\documentclass[11pt,leqno]{article}
\usepackage{amsfonts, amsthm}
\begin{document}
\newtheorem{exa}{Example}
\newtheorem{theo}{Theorem}
\newtheorem{rem}{Remark}
\newtheorem{lem}{Lemma}
\newtheorem{defi}{Definition}
\newtheorem{coro}{Corollary}
\newtheorem{prob}{P}
\def\proj{{\mathbb P}^2} 
\def\R{{\mathbb R}}
\def\Z{{\mathbb Z}}
\def\N{{\mathbb N}}
\def\C{{\mathbb C}}
\def\Q{{\mathbb Q}}
\newcommand{\tc}[2]{TC_{#1}{#2}}  
\newcommand{\ffc}[1]{{{\cal F}}(#1)}
\newcommand{\wc}[1]{{{\cal M}}(#1)}
\newcommand{\ps}[1]{{{\mathbb P}}^{#1}}
\def\F{{\cal F}}
\def\Nn{{\cal N}}
\def\P{{\cal P}}
\def\Pt{\C[t]}
\def\RR{{\cal R}} 
\def\li{{\mathbb L}}
\def\L{{\cal L}}
\newcommand{\cfc}[1]{\{\delta_t\}_{t\in #1}}

\begin{center}
{\LARGE\bf Center conditions for polynomial differential 
equations: discussion of some problems
\footnote{
\\
}
\\}
\vspace{.25in} {\large {\sc Hossein Movasati}} 
\\
Ochanomizu University \\
Department of Mathematics \\
2-1-1 Otsuka, Bunkyo-ku\\
Tokyo 112-8610, Japan \\
Email: movasati@cc.ocha.ac.jp
\end{center}
\begin{abstract}
Classifications of irreducible components of the set
of polynomial differential equations with a fixed degree
and with at least one center singularity lead to some
other new problems on Picard-Lefschetz theory and Brieskorn
modules of polynomials. In this article we explain these problems and
their connections to such classifications. 
\end{abstract}
\setcounter{section}{-1}
\section{Introduction}
The set of polynomial 1-forms  $\omega=P(x,y)dy-Q(x,y)dx, 
\deg P, \deg Q \leq d,\ d\geq 2$ is a vector space of finite dimension and we
denote by $\overline{\ffc{d}}$ its projectivization.
Its subset $\ffc{d}$ containing all $\omega$'s with 
$P$ and $Q$ relatively prime and 
$\deg(\omega):=\max\{\deg P,\deg Q\}=d$ is Zariski open
in $\overline{\ffc{d}}$. 
We denote the elements of $\overline{\ffc{d}}$ by $\F(\omega)$ or $\F$ if 
there is no confusion about the underlying 1-form $\omega$ in the text. 
Any $\F(\omega)$ induces a holomorphic foliation $\F$ in $\C^2$ i.e., 
the restrictions of $\omega$ to the leaves of $\F$ are 
identically zero. Therefore, we name an element of $\ffc{d}$ a
(holomorphic) foliation of degree $d$. 
 
The points in ${\mathrm sing} (\F(\omega))=\{P=0,Q=0\}$ are called the 
singularities  
of $\F(\omega)$. A singularity $p\in\C^2$ of $\F(\omega)$ is called reduced if 
$(P_xQ_y-P_yQ_x)(p)\not =0$. 
A  reduced singularity $p$ is called a center singularity or center for
simplicity if 
there is a holomorphic coordinates system $(\tilde x,\tilde y)$ around 
$p$ with  $\tilde x(p)=0,\tilde y(p)=0$ such that in this coordinates 
system $\omega\wedge d(\tilde x^2+\tilde y^2)=0$. One can call 
$f:=\tilde x^2+\tilde y^2$ a local first integral around $p$.  
The leaves of $\F$ around the center $p$ are given by 
$\tilde x^2+\tilde y^2=c$. 
Therefore, the leaf associated to the constant $c$ contains the one 
dimensional cycle 
$\{(\tilde x\sqrt{c},\tilde y\sqrt{c})\mid (\tilde x,\tilde y)\in \R^2,
\tilde x^2+\tilde y^2=1\}$ which is called
the vanishing cycle. We consider 
the subset of $\ffc{d}$ containing 
$\F(\omega)$'s with at least one center and we denote its closure
in $\overline{\ffc{d}}$ by $\wc{d}$.
It turns out that $\wc{d}$ is 
an algebraic subset of $\ffc{d}$ (see for instance \cite{mov0}). 
Now the problem of identifying the irreducible components of 
$\wc{d}$  arises. This problem is also known by the name 
"Center conditions" in the context of real polynomial differential 
equations. Let us introduce some of irreducible
components of $\wc d$.

For  $n\in\N\cup\{0\}$, 
let $\P_n$ denote the set of polynomials of degree at most 
$n$ in $x$ and
$y$ variables. 
Let also $d_i\in \N,\ i=1,2,\ldots,s$ with $\sum_{i=1}^s d_i=d-1$
and
$\L(d_1,\ldots,d_s)$ be the set of logarithmic foliations 
$$
\F(f_1\cdots f_s\sum_{i=1}^s \lambda_i\frac{df_i}{f_i}),\ f_i\in 
\P_{d_i},\ \lambda_i\in \C
$$
For practical purposes, one  assumes that 
$\deg f_i=d_i, \lambda_i\in\C^*,\  1\leq i \leq s$ and that 
$f_i$'s intersect each other transversally, and one  obtains an
element in $\ffc{d}$. Such a foliation has the logarithmic first 
integral $f_1^{\lambda_1}\cdots f_s^{\lambda_s}$. 
Since $\L(d_1,\ldots,d_s)$ is parameterized  by 
$\lambda_i$ and $f_i$'s it is irreducible.   
\begin{theo}(\cite{hosmov})
\label{main}
The set $\L(d_1,\ldots,d_s)$ is an irreducible component 
of $\wc{d}$, where $d=\sum_{i=1}^{s} d_i-1$.
\end{theo}
In the case $s=1$ we can assume that $\lambda_1=1$ and so $\L(d+1)$ is
the space of foliations of the type $\F(df)$, where $f$ is a polynomial 
of degree $d+1$. This case is proved by Ilyashenko in \cite{ily}.    

In general the aim is to find $d_i\in\N\cup\{0\}, i=1,2,\ldots,k$ 
and parameterize an irreducible component
$X=X(d_1,d_2,\ldots, d_k)$ of $\wc d$ by $\P_{d_1}
\times\P_{d_2}\times\cdots\times\P_{d_k}$. In the
above example $k=2s$ and $d_{s+1}=\cdots d_{2s}=0$.
 Once we have done 
this, we can reformulate the fact that $X$ is an irreducible component of
$\wc{d}$ in a meaningful way as follows:
\begin{theo}
\label{main1}
There exists an open dense subset $U$ of 
$X$ with the following property: 
for all $\F\in U$ parameterized with 
$f_i\in\P_{d_i},\ i=1,2,\ldots,k$ 
and a center $p\in\C^2$ of $\F$ let $\F_\epsilon$  be a holomorphic
deformation of $\F$ in $\ffc{d}$ such that its unique
singularity $p_\epsilon$ near $p$ is still a center.
Then there exist polynomials $f_{i\epsilon}\in \P_{d_i}$ 
such that $\F_\epsilon$ is parameterized by 
$f_{i\epsilon}$'s. Here $f_{i\epsilon}$'s are holomorphic in 
$\epsilon$ and  $f_{i0}=f_i$.
\end{theo}
The above theorem also says that the persistence of one center implies
the persistence of all other type of singularities.  
\section{Usual method}
To prove theorems like Theorem \ref{main1} usually 
one has to take $U$ the complement of 
$X\cap {\mathrm sing}(\wc{d})$ in $X$. But
this is not an explicite description of $U$. In practice
one defines $U$ by conditions like: $f_i,\ i=1,2,\ldots,k$ is of degree
$d_i$, $f_i$'s have no common factors, 
$\{f_i=0\}$'s intersect each other transversally and so on. To 
prove Theorem \ref{main1}, after
finding such an open set $U$, it is enough to prove that for at least one
$\F\in U$
\begin{equation}
\label{10.3.04}
T_\F X=T_\F\wc{d}
\end{equation}
where $T_\F$ means the tangent bundle at $\F$. 
Note that for a foliation $\F\in X$ the equality
(\ref{10.3.04}) does not imply that $\F\in U$. 
There may be an irreducible component of
$\wc{d}$ of dimension lower than the dimension of $X$ such that it passes 
through $\F$ and its tangent space at $\F$ is a subset of $T_\F X$. For
this reason  after proving (\ref{10.3.04}) for $\F$ with some
generic conditions on $f_i$'s,  we may not be sure that $U$ defined by 
such generic  conditions on $f_i$'s is $X-(X\cap {\mathrm sing}(\wc{d}))$. 
However, in the 
bellow $U$ can mean 
$X-(X\cap {\mathrm sing}(\wc{d}))$ or some open dense subset of $X$.  

An element $\F$ of the irreducible component $X$  may have more than one center.
The deformation of $\F$ within $X$ may destroy some centers but it preserves
at least one center. Therefore, we have the notion of stable and unstable center
for elements of $X$. A stable center of $\F$ is a center which persists 
after any  deformation of $\F$ within $X$. An unstable center is a center 
which is not stable. It is natural to ask
\begin{prob}
Are all the centers of a foliation $\F\in U$ stable?
\end{prob}
The answer is positive for $X=\L(d_1,d_2,\ldots,d_s)$ in Theorem \ref{main}.
Every element $\F\in U$ has $d^2-\sum_{i<j}d_id_j$ stable center. Here $U$
means just an open dense subset of $X$.  

The inclusion $\subset$ in the equality (\ref{10.3.04}) is trivial. 
To  prove the other side $\supset$, we fix a stable center singularity $p$ 
of $\F$ and make a deformation $\F_\epsilon(\omega+\epsilon\omega_1+\cdots)$ 
of  $\F=\F(\omega)$. Here $\omega_1$ represents an element $[\omega_1]$ of 
$T_\F\wc{d}$.  
Let $f$ be a local first integral in a neighborhood $U'$ of $p$, 
$s$ a holomorphic
function in $U'$ such that $\omega=s.df$,  
$\delta$ a vanishing cycle in a leaf of $\F$ in $U'$  and 
$\Sigma\simeq(\C,0)$ a transverse section to $\F$ in a 
point $p\in \delta$. We  assume that the transverse section
$\Sigma$ is parameterized by $t=f\mid_\Sigma$. 
The holonomy of $\F$ along $\delta$ is identity.
Let $h_\epsilon(t)$ be the holonomy of $\F_\epsilon$ along the 
path $\delta$. It is a holomorphic function in $\epsilon$ and $t$ and
 by hypothesis $h_0(t)=t$. We write the Taylor expansion of 
$h_\epsilon(t)$ in $\epsilon$
$$
h_\epsilon(t)-t=M_1(t)\epsilon+M_2(t)
\epsilon^2+\cdots +M_i(t)\epsilon^i+\cdots,\  i!.M_i(t)=
\frac{\partial^ih_\epsilon}{\partial \epsilon^i}\mid_{\epsilon=0}
$$
The function $M_i$ is called the $i$-th Melnikov function of the deformation
$\F_\epsilon$  along the path  $\delta$. It is well-known that 
the first Melnikov function is given by
$$
M_1(t)=-\int_{\delta_t}\frac{\omega_1}{s}
$$
where $\delta_t$ is the lifting up of 
$\delta$ in the 
leaf through $t\in \Sigma$, and the multiplicity 
of  $M_1$ at $t=0$ is the number of limit cycles 
(more precisely the number of fixed points
of the holonomy $h_\epsilon$) which appears around $\delta$ 
after the deformation  (see for instance \cite{mov0}). 
This fact shows the importance of these functions in the local study 
of Hilbert 16-th problem.

Now, if in the deformation $\F_\epsilon$ the deformed singularity 
$p_\epsilon$ near $p$ is center  then $h_\epsilon={\mathrm  id}$
and in particular 
\begin{equation}
\label{stuhler}
\int_{\delta_t}\frac{\omega_1}{s}=0,\ \forall t\in \Sigma
\end{equation}
Let $T^*_\F X$ be the set of $[\omega_1]\in T_\F\F(d)$ with
the above property. 
It is easy to check that the above definition does not depends on
the choice of $f$ (see \cite{mov0}). 
We have seen that $T_\F\wc d\subset T^*_\F X$. The following 
question arises:
\begin{prob}
\label{p1}
Is $T_\F\wc d=T^*_\F X$?
\end{prob}
If the answer is positive then it means that form the vanishing
of integrals (\ref{stuhler}) one must be able to prove that 
$\omega_1\in T_\F X$. Otherwise,
calculating more Melnikov functions to
get more and more information
on $\omega_1$ is necessary. The proof of Theorem \ref{main} with $s=1$ shows
that the answer of P\ref{p1} is positive in this case. However, the answer of
P\ref{p1} for
$X=\L(d_1,d_2,\ldots,d_s)$ is not known. 
\section{Some singularities of $\wc d$}
The method explained in the previous section
has two difficulties: First, 
identifying $U:=X\cap sing(\wc{d})$ 
and second to know the dynamics and topology of
the original foliation $\F$. A way to avoid these difficulties
is to look for foliations $\F(df)$, where $f$ is a degree
$d+1$ polynomial in $\C^2$. We already know that such foliations
lie in the irreducible component $\L(d+1)$. But if we take $f$ a 
non-generic polynomial then $\F(df)$ may lie in other irreducible components of $\wc d$ and
even worse, $\F(df)$ may not be a smooth point of such irreducible components.
\begin{prob}
\label{p2}
Do all irreducible components of $\wc{d}$ intersect $\L (d+1)$?
\end{prob}
If the answer of the above question is positive then the classification
of irreducible components of $\wc d$ leads to the classification
of polynomials of degree $d+1$ in $\C^2$ according to their Picard-Lefschetz
theory and Brieskorn modules. If not, we may be 
interested to find an irreducible
component $X$ which does not intersect $\L (d+1)$. In any case, the method
which we are going to explain bellow is useful for those $X$ which intersect
$\L(d+1)$. 

The foliation $\F=\F(df)$  has a first integral $f$ and so it has 
no dynamics. The function $f$ induces a ($C^\infty$) locally trivial 
fibration on $\C-C$, where $C$ is a  finite subset of $\C$. 
The points of $C$ are  called critical values of $f$ and the associated 
fibers are called the critical fibers. 
We have  Picard-Lefschetz theory of $f$ and the action of monodromy
$$
\pi_1(\C-C,b)\times H_1(f^{-1}(b),\Q)\rightarrow H_1(f^{-1}(b),\Q)
$$
where $b\in \C-C$ is a regular fiber. Let $\delta'\in H_1(f^{-1}(b),\Q)$ 
be the monodromy of $\delta$ (the vanishing cycle around a center singularity
of $\F(df)$) along an arbitrary path in $\C-C$ with the end point $b$.
From analytic continuation
of the integral (\ref{stuhler}) one concludes that 
$\int_{\pi_1(\C-C).\delta}\omega=0$.   
\begin{prob}
\label{pl}
Determine the subset $\pi_1(\C-C).\delta\subset H_1(f^{-1}(b),\Q)$.
\end{prob}
In the case of a generic polynomial $f$, Ilyashenko has proved that in 
P\ref{pl} 
the equality happens. To prove Theorem \ref{main}, I have used a 
polynomial $f$ which is a product of $d+1$ lines in general position and
I have proved that $\pi_1(\C-C).\delta$ together with the  cycles at infinity
generate $H_1(f^{-1}(b),\Q)$. Cycles at infinity are cycles around the points
of compactification of $f^{-1}(b)$.  

Parallel to the above topological theory theory, we have another algebraic
theory associated to each polynomial. The Brieskorn module 
$H=\frac{\Omega^1}{d\Omega^0+\Omega^0 df}$, where 
$\Omega^i,i=0,1,2$ is the set of polynomial differential $i$-forms in $\C^2$, 
is a $\C[t]$-module in a natural way and we have the action of Gauss-Manin
connection 
$$
\nabla : H_C\rightarrow H_C
$$
where $H_C$ is the localization of $H$ over the multiplicative
subgroup of $\C[t]$ generated by $t-c,\ c\in C$ (see \cite{hosmov}). 
\begin{prob}
\label{br}
Find the torsions of $H$ and classify the kernel of the maps
$\nabla^i=\nabla\circ\nabla\circ\cdots\circ\nabla $ $i$-times.
\end{prob}
When $f$ is the product of lines in general position then $H$ has
not torsions and the classification of the kernel of $\nabla^i$ 
is done in \cite{hosmov} using a theorem of Cerveau-Mattei.

Solutions to the both problems P\ref{pl} and P\ref{br} are 
closely related to the position of $\F(df)$ in $\wc{d}$.
Using solutions to P\ref{pl} and P\ref{br} one calculates the Melnikov 
functions $M_i$'s by means of 
integrals of 1-forms (the data of the deformation) over vanishing cycles 
and one calculates the tangent cone $TC_\F\wc{d}$ of $\F=\F(df)$ in
$\wc d$ and compare it with the tangent cone of suspicious irreducible
components of $\wc d$. For instance, to prove Theorem \ref{main},
we have taken $f$ the product of $d+1$ lines in general position
and we have proved that 
\begin{equation}
\label{tc}
\cup_{\sum_{i=1}^s d_i=d-1} TC_\F\L(d_1,d_2,\ldots,d_s)=TC_\F\wc d
\end{equation}
All the varieties $\L(d_1,\ldots,d_s),\ \sum_{i=1}^s d_i=d-1$ pass through 
$\F=\F(df)$.
\begin{prob}\rm
 Are $\L(d_1,\ldots,d_s)$'s  all irreducible components of $\wc{d}$
through $\F(df)$?
\end{prob}
 Note that the equality (\ref{tc}) 
does not give an answer to this problem. 
There may be an irreducible 
component of $\wc{d}$ through $\F(df)$ and different form
$\L(d_1,d_2,\ldots,d_s)$'s such that
its tangent cone at $\F(df)$ is a subset of (\ref{tc}). 
In this case the definition of other notions of
tangent cone based on higher order 1-forms in the deformation of 
$\F(df)$ seems to be necessary. 

The first case in which one may be interested to use the method of this section
can be:
\begin{prob}\rm
\label{2may02}
Let $l_i=0,\ i=0,1,\ldots, d$ be lines in the real plane and 
$m_i,\ i=0,1,\ldots, d$ be integer numbers. 
Put $f=l_0^{m_0}\cdots l_d^{m_d}$. Find all irreducible components 
of $\wc{d}$ through $\F(df)$. 
\end{prob}
In this problem the line $l_i$ has multiplicity $m_i$ and it would be
interesting to see how the classification of irreducible components through
$\F(df)$ depends on the different arrangements of the lines $l_i$ in the 
real plane and the associated multiplicities. In particular, we may allow
several lines to pass through a point or to be parallel. 
When there are lines with negative multiplicities then we have a third kind
of singularities $\{l_i=0\}\cap\{l_j=0\}$ called dicritical singularities, 
where $l_i$ (resp. $l_j$) has positive (resp. negative) multiplicity. 
They are indeterminacy points of $f$ and are characterized by this property 
that there are infinitely many leaves of the foliation passing through 
the singularity. Also in this case there are saddle critical points of $f$ 
which are not due to the intersection points of the lines with positive (resp.
negative) multiplicity. The reader may analyze the situation by the example
$f=\frac{l_0l_1}{l_2l_3}$.
\section{Looking for irreducible components of $\wc d$}
To apply the methods of previous sections one must find some
irreducible subsets of $\wc{d}$ and then one conjectures that they must
be irreducible components of $\wc{d}$. 
The objective of this section is to do this.

Classification of codimension one foliations
on complex manifolds of higher dimension is a subject related to center
conditions. We state the problem in
the case of $\C^n, \ n>2$ which is compatible with this text. However, 
the literature on this subject is mainly for projective spaces 
 of dimension greater than two (see \cite{celi}). 

The set of polynomial 1-forms  $\omega=\sum^n_{i=1}P_i(x)dx_i, 
\deg P_i\leq d$ is a vector space of finite dimension and we
denote by $\overline{\F(n,d)}$ its projectivization.
Its subset $\F(n,d)$ containing all $\omega$'s with 
$P_i's$ relatively prime and 
$\deg(\omega):=\max\{\deg P_i, i=1,2,\ldots, n\}=d$ is Zariski open
in $\overline{\F (n,d)}$. 
An element $[\omega]\in \overline{\F (n,d)}$ induces a holomorphic foliation 
$\F=\F(\omega)$ in $\C^n$ if and only if
$\omega$ satisfies the integrability condition
\begin{equation}
\label{inte}
\omega\wedge d\omega=0 
\end{equation}
This is an algebraic equation on the coefficients of $\omega$. Therefore,
the elements of $\F(n,d)$ which induce a holomorphic foliation in $\C^n$ form an
algebraic subset, namely $\wc {n,d}$, of $\F(n,d)$. Now we have the problem of
identifying the irreducible components of $\wc{n,d}$. We define 
$\F(2,d):=\F(d)$ and $\wc{2, d}:=\wc{d}$.

Let us be given a polynomial map $F:\C^2\rightarrow \C^n, \ n \geq 2$ and 
a codimension one foliation $\F=\F(\omega)$ in $\C^n$. In the case $n>2$, 
let us suppose that $F$ is regular in a point 
$p\in\C^2$. This implies that $F$ 
around $p$ is a smooth embedding. We assume that $F(\C^2,p)$ has
a tangency with the  leaf of $\F$ through $F(p)$. 
In the case $n=2$, we assume that $F$ is singular at $p$. 
In both cases, after choosing a generic $F$ and $\F$,  
the pullback of $\F$ by $F$ has a center singularity 
at $p$. 
\begin{prob}
Fix an irreducible component $X$ of $\F(n,d)$. Is $$
\{F^*\F,\F\in X,\ 
\deg f_i\leq d_i,\ i=1,2,\ldots, n\}
$$
where $F=(f_1,f_2,\ldots,f_n)$,
an irreducible component of $\wc {d''}$ for some $d''\in\N$? 
\end{prob}
For instance in Theorem \ref{main}, the elements of $\L(d_1,d_2,\ldots,d_s)$
are pull backs of holomorphic foliations $\F(x_1x_2\cdots x_s\sum_{i=1}^s
\lambda_i\frac{dx_i}{x_i}),\ \lambda_i\in\C^*$ in $\C^s$ by the polynomial
maps $F=(f_1,f_2,\ldots, f_s),\ \deg f_i\leq d_i$.

Another way to find irreducible subsets of $\wc {d}$ is by looking for
foliations of lower degree. Take a polynomial of degree $d$ in $\C^2$ 
with the generic conditions considered by
Ilyashenko, i.e. $f$ has non degenerated singularities with distinct images.
Now $\F(df)$ has degree $d-1$ which is less than the degree of a generic 
foliation in $\F(d)$. 
\begin{prob}\rm
Classify all irreducible components of $\wc{d}$ through $\F(df)$.
\end{prob}
All $\L(d_1,\ldots,d_s)$'s pass through $\F(df)$. 
There are other candidates as follows: 
\begin{enumerate}
\item
$A_i=\{\F(\frac{dp}{p}+d(\frac{q}{p^i}))\mid deg(p)=1,deg(q)=d\}$ $i=0,1,2,\ldots, d$;
\item
$B_1=\{\F(\frac{dq}{q}+d(p))\mid deg(p)=1,deg(q)=d\}$;
\end{enumerate}
An element of $A_i$ (resp. $B_1$) has a first integral of 
the type $pe^{q/p^i}$ (resp. $qe^p$). 
These candidates are supported by Dulac's classification (see \cite{dul} and
\cite{celi} p.601) in the case $d=2$.

We can look at our problem in a more general context. Let
$M$ be a projective complex manifold of dimension two. We consider
the space $\F(L)$ of holomorphic foliations in $M$ with the 
normal line bundle $L$ (see for instance \cite{mov0}). Let also
$\wc {L}$ be its subset containing holomorphic foliation with at least one
center singularity. Again $\wc L$ is an algebraic subset of $\F(L)$ and one can
ask for the classification of irreducible components of $\wc L$. 
For $M=\C P(2)$ some irreducible components of 
$\wc L$ are identified  in \cite{mov0}.
\begin{prob}\rm
Prove a theorem similar to Theorem \ref{main} for an arbitrary projective
manifold of dimension two.
\end{prob}
In this generality one must be careful about trivial centers which 
we explain now. Let
$\F$ be a holomorphic foliation in $\C^2$ and $0$ a regular point
of $\F$. We make a blow up (see \cite{casa}) 
at $0$ and we obtain a divisor ${\mathbb C}P(1)$
which contains exactly one singularity of the blow up foliation and
this singularity is a center. 


\end{document}